\documentclass[10pt]{article}
\usepackage{amsmath}
\usepackage{amssymb}

\allowdisplaybreaks

\numberwithin{equation}{section}

\def\eps{\varepsilon}

\def\qed{\hfill\rule{.2cm}{.2cm}}
\def\P{{\mathbb P}}
\def\E{{\mathbb E}}
\def\Z{{\mathbb Z}}

\def\p{{\cal P}}
\def\e{{\cal E}}

\def\w{{\mathfrak W}}

\def\a{{\alpha}}

\def\1{{\bf 1}}

\newtheorem{theo}{Theorem}[section]

\newtheorem{lm}{Lemma}[section]

\newtheorem{rmk}{Remark}[section]

\title{Repulsion of an evolving surface on walls with random heights}
\author{L.R.G. Fontes\thanks{Partially supported by CNPq
grants 300576/92-7 and 662177/96-7 (PRONEX) and FAPESP grant 99/11962-9.}$^{~,1}$
  \and M.~Vachkovskaia\thanks{Partially supported  FAPESP grant 99/11962-9 and CNPq
  grant 306029/2003-0.}$^{~,2}$
  \and A.~Yambartsev\thanks{Supported by FAPESP
   grant 02/01501-9.}$^{~,1}$}

\begin{document}
\maketitle
{\footnotesize

\noindent
$^1$Department of Statistics, Institute of Mathematics and
Statistics, University of S\~ao Paulo, rua do Mat\~ao 1010,
CEP 05508--090, S\~ao Paulo SP, Brazil.\\
E-mails: lrenato@ime.usp.br, yambar@ime.usp.br

\noindent
$^2$Department of Statistics,  Institute of Mathematics,
Statistics and Scientific Computation, University of Campinas,
Caixa Postal 6065, CEP 13083--970, Campinas SP, Brazil.\\
E-mail: marinav@ime.unicamp.br}

\begin{abstract}
We consider the motion of a discrete
random surface interacting by exclusion with a random wall. The heights of the wall
at the sites of $\Z^d$ are i.i.d.\ random variables. Fixed the wall configuration,
the dynamics is given by the serial
harness process which is not allowed to go below the wall.
We study the effect of the
distribution of the wall heights on the repulsion speed.



\end{abstract}

\section{Introduction}

This paper is part of a project aiming to understand the effect of
the interaction with walls on the evolution of a $d$-dimensional
random surface in $(d+1)$-space.

The evolving random surface is modeled by the harness process
introduced by Hammersley in~\cite{H}, where among other results,
the fluctuations of the free case (no wall) were established in
all dimensions (see also~\cite{ffnv}, where this is discussed in
more detail than in here).

In~\cite{ffnv}, a solid flat wall is placed at the origin and its
effect on the displacement of the surface with respect to its
initial location at the origin is studied.

In that reference, it is shown that in all dimensions the average
height of the surface (say, at the origin) diverges to $+\infty$
as time increases. This should be compared to the average absolute
height of the surface at the origin in the free case. In the
latter case, that quantity is bounded in dimensions 3 and
higher~\cite{H, ffnv}. An effect of repulsion on the wall is thus
established in those dimensions. Estimates on the speed of
repulsion are obtained for a class of noise distributions
(including the Gaussian case). These are comparable to estimates
of the entropic repulsion for the massless free field interacting
with a flat wall (see~\cite{ffnv} and references therein).

Motivated by work on the entropic repulsion for the massless free
field interacting with a wall with random heights~\cite{BG}, we
consider the same kind of wall here. In~\cite{BG}, estimates
similar to those in~\cite{BDG, BDZ, D} for the wall with fixed
height case were obtained, showing in some cases an effect of the
wall height distribution. We show the analogous effect, with
analogous quantitative estimates, for the same class of noise
distributions considered in~\cite{ffnv} (see Theorem~\ref{rheight}
below).

Further studies on the massless free field interacting with a wall
with random heights were carried on in~\cite{BG1,BG2}. We refer
again to~\cite{ffnv} for other works on surfaces interacting with
walls, in and out of equilibrium.

In the next section we define precisely our model and describe
the flat wall result of~\cite{ffnv}, which is related and relevant
to our main result. The latter is presented and argued in the following 
and final section. It was announced previously in~\cite{fvy}.


\section{The model}
\label{mod}

Denote by $|i-j|$ the  number of edges in a minimal path
connecting $i$ and $j$ (we will use this definition not only for
$\Z^d$, but  also for other graphs). Let ${\mathcal
P}=\{p(i,j)\}_{i,j \in \Z^d}$ be a symmetric stochastic matrix
which satisfies $p(i,j) = p(0, j-i)=:p(j-i)=p(i-j)$ (homogeneity)
and $p(j) = 0$ for all $|j| > v$ for some $v$ (finiteness). Assume
also that ${\mathcal P}$ is truly $d$-dimensional: $\{j\in
\Z^d\,:\, p(j)\neq 0\}$ generates $\Z^d$. The weights  $p(i,j)$
can be interpreted  as transition probabilities of a random walk
on $\Z^d$; denote by  $\p$ its transition matrix and by $p_m(i,j)$
its $m$-step transition probabilities. By homogeneity,
$p_m(i,j)=p_m(0,j-i)=:p_m(j-i)$.

Let $\e:=\{\eps, \eps_n(i),\, i \in \Z^d, n \in \Z\}$ be a family
of i.i.d.~integrable symmetric random variables with unbounded
support. $\e$ represents the evolution noise variables.

We next introduce the wall variables, giving the heights at each
space coordinate. Consider the family of i.i.d.\ random variables
${\mathfrak W}=\{W(i)\}_{i\in\Z^d}$, independent of $\e$. $W(i)$
represents the height of the wall at site $i$.

With a realization of ${\mathfrak W}$ fixed, the harness process
interacting with ${\mathfrak W}$ by exclusion is defined as
follows.
\begin{equation}
\label{W2} X^\w_n(i) =
\begin{cases}
{}\hspace{2cm}0\vee W(i), & \mbox{ if } n =0,\\
W(i)+\Big(\p X^\w_{n-1}(i) +\eps_n(i)-W(i)\Big)^+, &\mbox{ if
}n\geq1.
\end{cases}
\end{equation}
We allow $W(i)$ to take the value $-\infty$ (with positive
probability), in which case the expression for $n\geq1$
in~(\ref{W2}) is
\begin{equation}
  \label{eq:W2l}
  \p X^\w_{n-1}(i) +\eps_n(i).
\end{equation}

\begin{rmk}
\label{ord} Notice that $X^\w_n(i)$ so defined is nondecreasing in
(the natural partial order for) $\w$.
\end{rmk}

The case where ${\mathfrak W}\equiv-\infty$, in which we denote
$X^\w$ by $Y$, is the free case introduced by Hammersley
in~\cite{H}. In that paper it is shown that $\E Y_n^2(0)$ is of
order $1$ in $d\geq3$. (Notice that under our assumptions on
$X^\w_0$ and $\eps$, $\E Y_n(0)\equiv0$.)

The case where $\w\equiv0$, in which we denote $X^\w$ by $Z$, was
studied in~\cite{ffnv}. We now quote some of the results of that
paper, which are directly related to our main result here.
\begin{theo}
\label{bounds_gauss} {\it (Part of Theorem 1.2
from~\cite{ffnv}).}\/ Let $F_\eps$ be the distribution function of
$\eps$ and define the following classes of distribution functions:
\begin{eqnarray}
  {\mathcal L}_\a^-&:=&\{F:\bar F(x)\leq ce^{-c'x^{\a}}, x>0,
                     \mbox{ for some positive $c$, $c'$}\},\label{eq:logpol-}\\
 {\mathcal L}_\a^+&:=&\{F:\bar F(x)\geq ce^{-c'x^{\a}}, x>0,
                     \mbox{ for some positive $c$, $c'$}\},\label{eq:logpol+}
\end{eqnarray}
where $\bar F=1-F$, and
\begin{eqnarray}
{\mathcal L}_\a&:=&{\mathcal L}_\a^-\cap{\mathcal L}_\a^+.
\label{eq:logpol}
\end{eqnarray}
For $d\geq3$, there exist constants $c$ and $C$ that may depend on
the dimension such that
\begin{itemize}
\item[(i)] if $F_\eps\in {\mathcal L}_\a$ for some
$1\le\a\ne1+d/2$, then
\begin{equation}
    \label{eq:b3}
    c (\log n)^{\frac1\a}\leq  \E\; Z_n(0)  \leq C(\log n)^{\frac1\a\vee\frac2{2+d}};
  \end{equation}
\item[(ii)] if $F_\eps\in {\mathcal L}_{1+d/2}$, then
\begin{equation}
    \label{eq:b4}
   c (\log n)^{\frac{2}{2+d}}\leq\E\; Z_n(0)\leq C (\log
   n)^{\frac{2}{2+d}}(\log\log n)^{\frac{d}{2+d}}.
  \end{equation}
\end{itemize}
\end{theo}


\section{Results}
\label{randomheight}

In the following, which is our main result, we obtain bounds on
the average height of the wall at the origin as a function of $n$,
the number of iterations of the dynamics, in $d\geq3$. The average
is taken with respect to the noise and wall variables. The bounds
are similar to the corresponding ones in
Theorem~\ref{bounds_gauss} above, and show an effect of the wall
variables (to leading order, ignoring constants) when they have a
heavy enough positive tail which is heavier than the noise ones.
This is the case when the noise variables are Gaussian and the
wall ones are sub-Gaussian (i.e., have distribution function
belonging to ${\cal L}_\theta$ with $\theta<2$).

\begin{theo}
\label{rheight}
Let $d\ge 3$ and
suppose that $F_\eps(x)\in {\cal L}_\alpha$
and $F_W(x)\in{\cal L}_\theta$.
 Then there exist $c$ and $C$ such that
\begin{equation}
\label{rhthm}
c(\log n)^{\frac{1}{\alpha}\vee \frac{1}{\theta}}\le \E\;X^\w_n(0)
\le C (A_n+(\log n)^{\frac{1}{\theta}}),
\end{equation}
where
\begin{eqnarray*}
A_n  = \left\{ \begin{array}{ll}
 (\log n)^{\frac{1}{\a}\vee\frac{2}{2+d}},&
\mbox{ if } \alpha\ne 1+\frac{d}{2}\\
 (\log n)^{\frac{2}{2+d}} (\log \log n)^{\frac{d}{d+2}},&
\mbox{ if }\alpha= 1+\frac{d}{2}.
\end{array}
\right.
\end{eqnarray*}
The lower bound in~(\ref{rhthm}) is valid for $\alpha, \; \theta>0$,
and the upper bounds are valid for
$\alpha>1$, $\theta>0$.
\end{theo}

\medskip
\noindent
{\it Proof.}

\noindent
{\it Lower bound.}\/
 Let next $\hat\w =\{\hat W(i)\}_{i\in\Z^d}$, where
 \[
\hat W(i)=\left\{\begin{array}{rcc}
-\infty, & \mbox{ if } & W(i)<0,\\
0, & \mbox{ if } & W(i)\ge 0.
\end{array}
\right.
\]
It then follows that $X^\w_n\geq X^{\hat\w}_n$.
So, we need to obtain a lower bound for $\mu_n:=\E\;X^{\hat\w}_n(0)$.

\begin{lm}
\label{const_dens}
We have $\mu_n\ge c\,(\log n)^{1/\alpha}$, where $c$ is a positive constant.
\end{lm}

\noindent
{\it Proof.}\/
Denote $q=\P(W(i)\geq 0)$. With a slight abuse of notation we
identify below $\hat\w$ with the set
$\{i\in\Z^d:\; \hat W(i)=0\}$.
We have then
\begin{eqnarray}
\mu_n&=&\E (X^{\hat\w}_n(0)\mid 0\in {\hat\w})q+\E (X^{\hat\w}_n(0)\mid 0\notin {\hat\w}) (1-q)\nonumber\\
&=&\E [\E( ({\cal P} X^{\hat\w}_{n-1}(0)+\eps_n(0))^+\mid {\hat\w}, \; 0\in {\hat\w})]q
     \nonumber\\&&
\,+\,\E [\E( {\cal P} X^{\hat\w}_{n-1}(0)+\eps_n(0)\mid {\hat\w}, \; 0\notin {\hat\w})](1-q)\nonumber\\
&=& \E {\cal P} X^{\hat\w}_{n-1}(0)+\E[\E(({\cal P} X^{\hat\w}_{n-1}(0)+\eps_n(0))^-\mid {\hat\w},0\in {\hat\w})]
                                          q\nonumber\\
&=&  {\cal P} \E X^{\hat\w}_{n-1}(0)+ \E[\E((-{\cal P} X^{\hat\w}_{n-1}(0)+\eps_n(0))^+\mid {\hat\w},0\in {\hat\w})]
                            q\nonumber\\
&\ge& \label{Jens} {\cal P} \E X^{\hat\w}_{n-1}(0)+G({\cal P} \E X^{{\hat\w}^0}_{n-1}(0))q\\
&=& \label{osn_mu} \mu_{n-1}+G({\cal P} \E X^{{\hat\w}^0}_{n-1}(0))q,
\end{eqnarray}
where $a^-:=a^+ -a$,   $G(x)=\E (\eps-x)^+ $, ${\hat\w}^0={\hat\w}\cup\{0\}$,
and (\ref{Jens}) is due to Jensen's inequality.

We want now to estimate $\E X^{{\hat\w}^0}_n(j)$ in terms of $\E X^{{\hat\w}}_n(j)$.
Consider the processes $X^{{\hat\w}^0}_n$, $X^{\hat\w}_n$ and $Y_n=X^\emptyset_n$
(free process), all coupled together by using the same $\eps_k(j)$.
We have  that $X^{{\hat\w}^0}_n(j)\ge X^{\hat\w}_n(j)\ge Y_n(j)$ for all $j$.
So, if $j\ne 0$, using~(\ref{W2}), (\ref{eq:W2l}) and the fact that for $a\ge c$ it holds
$(a+b)^+ -(c+b)^+\le a-c$,
 we get
\begin{equation}
\label{notW}
X^{{\hat\w}^0}_n(j)- X^{\hat\w}_n(j)\le\sum_{k\in\Z^d}p(j,k)(X^{{\hat\w}^0}_{n-1}(k)- X^{\hat\w}_{n-1}(k)).
\end{equation}
For $j=0$ we have
\begin{equation}
\label{Wi}
X^{{\hat\w}^0}_n(0)- X^{\hat\w}_n(0)\le\sum_{k\in\Z^d}p(0,k)(X^{{\hat\w}^0}_{n-1}(k)-
X^{\hat\w}_{n-1}(k))+\eps_n(0)^- +{\cal P}Y_{n-1}^-(0)
\end{equation}
(here we used the fact that for $a\ge c\ge g$ it holds $(a+b)^+-(c+b)\le a-c+b^-+g^-$).
Iterating~(\ref{notW}) and~(\ref{Wi}), one can get that
\begin{equation}
\label{notW-iter}
X^{{\hat\w}^0}_n(j)- X^{\hat\w}_n(j)\le\sum_{m=1}^n p_{m-1}(j)(\eps_{n-m+1}(0)^- +{\cal P}Y_{n-m}^-(0) ).
\end{equation}
Note that, as $d\ge 3$, random walk with transition matrix ${\cal P}$ is transient and
also $\E Y_{n}^-(0)\le const$, so~(\ref{notW-iter}) implies that
\begin{equation}
\label{mu_i}
\E X^{{\hat\w}^0}_n(j)-\E X^{\hat\w}_n(j)\le C_0
\end{equation}
for some $C_0>0$,
for all $n$, $i$, and $j$.

So, $\mu_n$ satisfies
\begin{equation}
\label{mu**}
\mu_n\ge\mu_{n-1}+{\tilde G}( \mu_{n-1})
\end{equation}
where ${\tilde G}(x)=G(x+C_0)q$.
A lower bound of $O((\log n)^{1/\alpha})$
for $\E X^\w_n(0)$ then follows as in the proof of Theorem~3.1 and Corollary~3.3
in~\cite{ffnv}.
\qed

\medskip
We derive next a lower bound of $O((\log n)^{1/\theta})$.
Let $\mu^\w_n(i)=\E(X^\w_n(i)|\w)$.
As  the dynamics  of the process can be re-written
as
\begin{eqnarray}
\label{W3}
X^\w_n(i) =\begin{cases}{}\hspace{2cm} 0\vee W(i), & \mbox{ if } n =0, \\
            {\cal P}X^\w_{n-1}(i) +  \eps_n(i)
                + \Big( W(i)-{\cal P} X^\w_{n-1}(i) -  \eps_n(i)\Big)^+,
                        & \mbox{ if } n \geq 1, \\
                      \end{cases}
\end{eqnarray}
we have
\begin{eqnarray}
\mu^\w_n(i)&=&{\cal P}\mu^\w_{n-1}(i)+
 \E(( W(i)-\sum_{j \in \Z^d} p(i,j) X^\w_{n-1}(j) -  \eps_n(i))^+\mid \w)\nonumber\\
 &\ge& \label{muW}
 {\cal P}\mu^\w_{n-1}(i)+(W(i)-{\cal P}\mu^\w_{n-1}(i)   )^+.
\end{eqnarray}

 For $\w$ fixed, and $i\in\Z^d$ let
\begin{equation}\label{nuW}
\nu^{\w}_n(i)=
\begin{cases}{}\hspace{2cm}
 0\vee W(i), & \mbox{ if } n =0,\\
 {\cal P}\nu^{\w}_{n-1}(i)+(W(i)-{\cal P}\nu^{\w}_{n-1}{(i)})^+,
      & \mbox{ if } n \geq 1.
\end{cases}
\end{equation}
 We then have (since $x+(a-x)^+$ is nondecreasing in $x$ for all $a$)
that $\mu^\w_n(i)\geq\nu^{\w}_n(i)$ for all $\w,n,i$.

 Let next $\tilde \w =\{\tilde W(i)\}_{i\in\Z^d}$, where
\begin{eqnarray*}
\tilde W(j)=\left\{ \begin{array}{cl}
W(j), &\mbox{ if } W(j)\geq0,\\
-\infty,  &\mbox{ if } W(j)<0.
\end{array}
\right.
\end{eqnarray*}
 It follows that $\nu^\w_n(i)\geq\nu^{\tilde\w}_n(i)$ for all $\w,n,i$.

 We will estimate $\nu^{\tilde\w}_n(i)$.
Let us decompose $\tilde\w$ in the following way.
$\tilde\w=\tilde\w^i\vee\tilde\w_i$, with
$\tilde\w_i=\{\tilde W_i(j)\}_{j\in\Z^d}$ and
$\tilde\w^i=\{\tilde W^i(j)\}_{j\in\Z^d}$,
where
\begin{equation}
\tilde W^i(j)=\begin{cases}
\tilde W(i), &\mbox{ if } i\ne j,\\
-\infty,  &\mbox{ if } i=j,
\end{cases}
\quad\tilde W_i(j)=\begin{cases}
-\infty, &\mbox{ if } i\ne j,\\
\tilde W(i),  &\mbox{ if } i=j.
\end{cases}
\end{equation}

\begin{lm}
For all $n$ and $j$ it holds
\begin{equation}
\label{nuWi}
\nu^{\tilde\w^i}_n(j)\vee \nu^{\tilde\w_i}_n(j) \le \nu^{\tilde\w}_n(j)
 \le \nu^{\tilde\w^i}_n(j)+\nu^{\tilde\w_i}_n(j).
\end{equation}
\end{lm}

\medskip
\noindent
 {\it Proof.}\/
We prove the lemma by induction. For $n=0$~(\ref{nuWi}) is evident.

Suppose~(\ref{nuWi}) holds for $n-1$.
For $j\ne i$, we have
\[\nu^{\tilde\w}_n(j)= {\cal P}\nu^{\tilde\w}_{n-1}(j)+
(\tilde W(j)-{\cal P}\nu^{\tilde\w}_{n-1}(j)   )^+,\]
\[\nu^{\tilde\w^i}_n(j)= {\cal P}\nu^{\tilde\w^i}_{n-1}(j)+
(\tilde W(j)-{\cal P}\nu^{\tilde\w^i}_{n-1}(j))^+,\]
and
\[\nu^{\tilde\w_i}_n(j)= {\cal P}\nu^{\tilde\w_i}_{n-1}(j).\]
Note that induction assumption implies ${\cal P}\nu^{\tilde\w^i}_{n-1}(j)\le
{\cal P}\nu^{\tilde\w}_{n-1}(j)$ and ${\cal P}\nu^{\tilde\w_i}_{n-1}(j)\le
{\cal P}\nu^{\tilde\w}_{n-1}(j)$. So,
 $\nu^{\tilde\w_i}_n(j) \le \nu^{\tilde\w}_n(j)$, and,  as
$x+(a-x)^+$ is increasing,  $\nu^{\tilde\w^i}_n(j) \le \nu^{\tilde\w}_n(j)$.
Also by induction assumption,
\[{\cal P}\nu^{\tilde\w}_{n-1}(j)\le
{\cal P}\nu^{\tilde\w^i}_{n-1}(j)+ {\cal
P}\nu^{\tilde\w_i}_{n-1}(j).\] As  $(a-x)^+$ is decreasing, we have
\[(\tilde W(j)-{\cal P}\nu^{\tilde\w}_{n-1}(j)   )^+\le
 (\tilde W(j)-{\cal P}\nu^{\tilde\w^i}_{n-1}(j))^+ .\]
 Thus,
 \[\nu^{\tilde\w}_n(j)
 \le \nu^{\tilde\w^i}_n(j)+\nu^{\tilde\w_i}_n(j).\]
The case $j=i$ is
similar. \qed

\medskip

Let us now estimate
${\cal P}\nu^{\tilde\w_0}_n(0)$. We suppose that $W(0)=:W>0$,
otherwise $\nu^{\tilde\w_0}_n\equiv0$. We have
\begin{equation}\label{nuWz}
\nu^{\tilde\w_0}_n(i)=
\begin{cases}
{}\hspace{1.6cm} W, & \mbox{ if } i=0,\, n =0,\\
{}\hspace{1.7cm} 0, & \mbox{ if } i\ne 0,\, n =0,\\
{\cal P}\nu^{\tilde\w_0}_{n-1}(0)+(W-{\cal P}\nu^{\tilde\w_0}_{n-1}{(0)})^+,
      & \mbox{ if } i=0,\, n \geq 1,\\
{}\hspace{1.35cm} {\cal P}\nu^{\tilde\w_0}_{n-1}(i),
      & \mbox{ if } i\ne 0,\, n \geq 1.
\end{cases}
\end{equation}

It follows readily by induction that
$\nu^{\tilde\w_0}_n(i)\leq W$ for all $i,n$. Thence, we have
$\nu^{\tilde\w_0}_n(0)\equiv W$. It is now readily verified by induction that
for $i\ne0$
\begin{equation}
  \label{eq:expr}
  \nu^{\tilde\w_0}_n(i)=W\sum_{k=1}^np_k^{\{0\}}(i,0),
\end{equation}
where, for any $j$, $p_k^{\{0\}}(j,0)$ is the probability that the
random walk with transition matrix
${\cal P}$ starting from $j$ returns to $0$ for the first time at step $k$.
It follows that
\begin{equation}
\label{PnuWn}
{\cal P}\nu^{\tilde\w_0}_n(0)=W\sum_{k=1}^{n+1}p_k^{\{0\}}(0,0)\le a W,
\end{equation}
where $a=\sum_{k=1}^\infty p_k^{\{0\}}(0,0)<1$, since $d\ge 3$.

By~(\ref{nuWi}) and~(\ref{PnuWn}), we have
\begin{equation}
\label{indep}
\nu^{\tilde\w}_n(0)\ge {\cal P}\nu^{\tilde\w}_{n-1}(0)+
((1-a)W-{\cal P}\nu^{{\tilde\w'}}_{n-1}(0))^+,
\end{equation}
where $\tilde\w'=\{\tilde W'(i)\}_{i\in\Z^d}$, with
$\tilde W'(i)=\tilde W(i)$ if $i\ne0$, and $\tilde W'(0)$ independent
of $\tilde\w$.
Taking now expectations with respect to $\tilde\w,\tilde W'(0)$,
and applying Jensen's inequality, we get
\begin{eqnarray}\nonumber
\E\nu^{\tilde\w}_n(0)\!\!\!&\ge&\!\!\!\E{\cal P}\nu^{\tilde\w}_{n-1}(0)
+{\mathbb G}(\E{\cal P}\nu^{\tilde\w'}_{n-1}(0))
={\cal P}\E\nu^{\tilde\w}_{n-1}(0)
+{\mathbb G}({\cal P}\E\nu^{\tilde\w'}_{n-1}(0))\\
\label{Enu}
&=&\!\!\!{\cal P}\E\nu^{\tilde\w}_{n-1}(0)
+{\mathbb G}({\cal P}\E\nu^{\tilde\w}_{n-1}(0))=\E\nu^{\tilde\w}_{n-1}(0)
+{\mathbb G}(\E\nu^{\tilde\w}_{n-1}(0)),
\end{eqnarray}
where ${\mathbb G}(x)=\E((1-a)W-x)^+$, and we have used the
equidistribution of $\tilde\w$ and $\tilde\w'$, and the translation
invariance of the joint distribution of $\tilde\w$ and $\e$.

We then see that $\nu_n\equiv\E\nu^{\tilde\w}_{n}(0)$ is of the same form
as $(3.3)$ in~\cite{ffnv}, and a lower bound of $O((\log n)^{1/\theta})$
for $\nu_n$ follows as in the proof of $(3.4)$ in~\cite{ffnv}.

\bigskip

\noindent {\it Upper bound.}\/ As in~\cite{ffnv}, in order to
obtain an upper bound, we compare the wall process with a free
process started sufficiently high. Let $X^{\w, r_n}_n$ and
$Y^{r_n}_n$ have the same evolution as $X^\w_n$ (resp. $Y_n$), but
$X^{\w, r_n}_0(i)=r_n\vee W(i)$, $Y^{r_n}_0(i)\equiv r_n$. Let
${\cal R}_n=\max\{W(i): \; |i|\le v n\}$ (recall that $v$ is a
range of ${\cal P}$),  and
\begin{eqnarray}
\label{alt}
a_n=\left\{ \begin{array}{ll}
2K ((\log n)^{\frac{1}{\a}\vee\frac{2}{2+d}}+(\log n)^{\frac{1}{\theta}}),&
\mbox{ if } \alpha\ne 1+\frac{d}{2}\\
2 K  ((\log n)^{\frac{2}{2+d}} (\log \log n)^{\frac{d}{d+2}}+(\log n)^{\frac{1}{\theta}}) ,&
\mbox{ if }\alpha= 1+\frac{d}{2}.
\end{array}
\right.
\end{eqnarray}
Note that
$\P({\cal R}_n>K(\log n)^{\frac{1}{\theta}})\le n^{c_1-c_2 K^{\theta}}$.
Take $r_n=a_n/2$.
We have
\begin{eqnarray}
\P(X^\w_n(0)\ge a_n)&\le&\P(X^{\w, r_n}_n(0)\ge a_n)\nonumber\\
&=&\P(X^{\w, r_n}_n(0)\ge a_n, X^{\w, r_n}_n(0)= Y^{r_n}_n(0) )\nonumber\\
&&\!\!\!\!\!
+\P(X^{\w, r_n}_n(0)\ge a_n, X^{\w, r_n}_n(0)\ne Y^{r_n}_n(0))\nonumber\\
&\le&
\label{UB}
 \P(Y^{r_n}_n(0)\ge a_n)+ \P( X^{\w, r_n}_n(0)\ne Y^{r_n}_n(0)).
\end{eqnarray}
In Section 5 of~\cite{ffnv} it was shown that
$\P(Y^{r_n}_n(0)\ge a_n)\le k n^{c_3-c_4 K}$.
As for $\P( X^{\w, r_n}_n(0)\ne Y^{r_n}_n(0))$, note that
$ X^{\w, r_n}_n(0)$ and $ Y^{r_n}_n(0))$ can be different if
either if ${\cal R}_n>r_n$, or if it occurs
\[\{Y^{r_n}_l(j)< {\cal R}_n\mbox{ for some }(l,j)\mbox{ with }
                l\leq n,\; |j|\leq v(n-l)\}.\]
We have then
\begin{eqnarray}
\P(Y^{r_n}_l(j)< {\cal R}_n)&\le& \P(Y^{r_n}_l(j)<K (\log n)^{\frac{1}{\theta}} )+
    \P({\cal R}_n> K (\log n)^{\frac{1}{\theta}})\nonumber\\
&\le&\P(Y_l(j)<K (\log n)^{\frac{1}{\theta}}-r_n)+n^{c_1-c_2 K^\theta}\nonumber\\
&\le &\P(Y_l(j)>r_n-K (\log n)^{\frac{1}{\theta}})+n^{c_1-c_2 K^\theta}\nonumber\\
&\le& kn^{c_5-c_6 K^{\theta\wedge 1}},
\end{eqnarray}
where  $\theta\wedge 1=\min\{\theta,1\}$, and
\[\P( X^{\w, r_n}_n(0)\ne Y^{r_n}_n(0))\le n^{c_1-c_2 K^\theta}+
\sum_{l=0}^n\sum_{|j|\le v(n-l)} kn^{c_5-c_6 K^{\theta\wedge 1}}
\le k' n^{c_7-c_8 K^{\theta\wedge 1}}.\]
So,
\[
\P(X^\w_n(0)\ge a_n)\le k^*n^{c'-c'' K^{\theta\wedge 1}}
\]
and, by taking $K$ large enough, this implies that
$\E X^\w_n(0)\le \frac{C}{2K} a_n$
(see end of Section 6 of~\cite{ffnv} for the reasoning in a similar
situation).
\qed

\end{document}